%% file: MPConvex7.tex
\documentclass{amsart}

\usepackage{amsfonts,amssymb,mathrsfs}
\usepackage{graphics,color,hyperref}

\newtheorem{theorem}{Theorem}[section]
\newtheorem{definition}[theorem]{Definition}
\newtheorem{corollary}[theorem]{Corollary}
\newtheorem{proposition}[theorem]{Proposition} 
\newtheorem{propdef}[theorem]{Proposition-Definition} 
\newtheorem{lemma}[theorem]{Lemma} 
\theoremstyle{remark}
\newtheorem{remark}[theorem]{Remark} 
\newtheorem{example}[theorem]{Example}
\title[The Minkowski Theorem for Max-plus Convex Sets]{The Minkowski Theorem for Max-plus Convex Sets}
\author{St\'ephane Gaubert}
\address{INRIA, Domaine de Voluceau, 78153, Le Chesnay C\'edex, France. T\`el: +33 1 39 63 52 58, Fax: +33 1 39 63 57 86}
\email{Stephane.Gaubert@inria.fr}
\author{Ricardo D. Katz}
\address{CONICET. Postal address: Dep. of Mathematics, Universidad Nacional de Rosario, Avenida Pellegrini 250, 2000 Rosario, Argentina.}
\email{rkatz@fceia.unr.edu.ar}
\date{May 2, 2006}

\DeclareMathAlphabet{\mathbbold}{U}{bbold}{m}{n}

\newcommand{\zero}{\mathbbold{0}}
\newcommand{\unit}{\mathbbold{1}}
\newcommand{\set}[2]{\left\{#1\mid\,#2\right\}}
\newcommand{\R}{\mathbb{R}}
\newcommand{\N}{\mathbb{N}}
\newcommand{\Rm}{\R\cup\{-\infty\}}
\newcommand{\rmax}{\R_{\max}}
\newcommand{\co}{\mbox{\rm co}\,}
\newcommand{\cone}{\mbox{\rm cone}\,}
\newcommand{\Rec}{\mbox{\rm rec}\,}
\newcommand{\supp}{\mbox{\rm supp}\,}
\newcommand{\ext}{\mbox{\rm ext}\,}
\newcommand{\Gcone}[1]{C_{#1}}
\newcommand{\clo}{\mbox{\rm clo}\,}
\newcommand{\eg}{\mbox{\rm ext-g}\,}

\begin{document}
\begin{abstract}
We establish the following max-plus analogue of Minkowski's theorem.
Any point of a compact max-plus convex subset of $(\Rm)^n$ 
can be written as the max-plus convex combination
of at most $n+1$ of the extreme points of this subset.
We establish related results for closed max-plus convex cones
and closed unbounded max-plus convex sets.
In particular, we show that a closed max-plus convex set can be
decomposed as a max-plus sum of its recession cone and
of the max-plus convex hull of its extreme points.
\end{abstract}

\maketitle

\section{Introduction}
The {\em max-plus segment} joining two points $u,v\in (\Rm)^n$
is the set of vectors of the form
$(\alpha + u)\vee (\beta +v)$ 
where $\alpha$ and $\beta$ are elements of $\Rm$ such
that $\alpha \vee \beta=0$.
Here, $\vee$ denotes the maximum of scalars, or the pointwise maximum
of vectors,
and for all scalars $\alpha\in \Rm$ and vectors $u\in (\Rm)^n$,
$\alpha+ u$ denotes the vector with entries $\alpha+u_i$.

A subset of $(\Rm)^n$ is {\em max-plus convex} if it contains
any max-plus segment joining two of its points. 
The {\em max-plus convex cone} generated by $u,v$ is the set of vectors
of the form $(\alpha + u)\vee (\beta +v)$, where $\alpha$ and $\beta$ are
arbitrary elements of $\Rm$. A subset of $(\Rm)^n$ is a
{\em max-plus convex cone} if it contains any max-plus convex
cone generated by two of its points. These definitions
are natural if one considers the max-plus semiring, which
is the set $\Rm$ equipped with the addition
$(a,b)\mapsto a\vee b$ and the multiplication
$(a,b)\mapsto a+b$. Max-plus convex cones are
also called {\em semimodules} over the max-plus semiring.
An example of max-plus convex set is given in Figure~\ref{figure1}:
the convex set $A$ is the closed grey region, together with the
portion of vertical line joining the point $b$ to it. Three max-plus
segments in general position, joining the pairs of points $(f,g)$,
$(h,i)$, and $(j,k)$, are represented in bold. By comparing
the shapes of these segments with the shape of $A$, one can
check geometrically that $A$ is convex.

In this paper we give representation theorems, in terms of 
extreme points and extreme rays, for max-plus convex sets 
and cones. 

Motivations to study the max-plus analogues of convex cones and convex 
sets arise from several fields, let us review some of these
motivations.
 
Max-plus convex sets were introduced by 
K. Zimmermann~\cite{zimmermann77}. Convexity is a powerful 
tool in optimization, and so, max-plus convex sets
arose in the quest of solvable optimization 
problems~\cite{zimmermann84,zimmermann03}.
See also the book of U. Zimmermann~\cite{zimmermann81}
for an overview.

Max-plus convex cones have been 
studied in idempotent analysis, after the 
observation due to Maslov that the solutions of 
an Hamilton-Jacobi equation associated with a deterministic
optimal control problem satisfy a ``max-plus'' 
superposition principle, and so, belong to 
structures similar to convex cones, which
are called semimodules or idempotent linear spaces~\cite{litvinov00,cgq02}. 
Such structures have been used, for instance, to characterize
the sets of stationary solutions of deterministic
optimal control problems~\cite{agw04b},
and to design numerical algorithms~\cite{mceneaney,agl},
to mention a recent application.

Max-plus convex cones have also been studied in relation 
to discrete event systems. The reader may consult
the survey papers~\cite{maxplus97,ccggq99} for
more background. In particular, reachable 
and observable spaces of certain timed discrete event 
systems are naturally equipped with structures of 
max-plus polyhedral cones~\cite{katz05}.
Earlier discrete event systems motivations have been at the
origin of the works~\cite{CGQ96a,CGQ97a,gaubert98n},
in which the theory of max-plus polyhedral cones has
been developed.

Of course, another interest in max-plus convexity stems
from abstract convex analysis~\cite{ACA}. Several recent
papers in this field, in particular
those of Mart{\'{\i}}nez-Legaz, Rubinov, and Singer~\cite{singer}, 
and Akian and Singer~\cite{akiansinger}, are related to max-plus
algebra.
   
A renewed interest in max-plus convex cones, or 
``tropical convex sets'', and specially, in tropical
polyhedra, has recently arisen 
in relation to tropical geometry
(in this context, ``tropical'' is essentially 
used as a synonym of ``max-plus'', or rather, of
the dual term, ``min-plus'').
Tropical analogues of polytopes
have been considered by Develin and Sturmfels~\cite{sturmfels03}, 
and also by 
Joswig~\cite{joswig04}
(the tropical polytopes they consider are special
finitely generated max-plus convex
cones, in which the generators have finite entries).
Develin and Sturmfels have also pointed out an
elegant relation between tropical polytopes and phylogenetic 
analysis~\cite{sturmfels03}.

Some of these motivations have guided the development
of max-plus analogues of classical results of convex
analysis, like the Hahn-Banach theorem~\cite{zimmermann77,shpiz,cgq02,cgqs03}. 

We are interested here in the representation
of convex sets in terms of extreme points or extreme rays.
This problem, in the case of finitely generated max-plus convex cones,
has been considered by several authors~\cite{moller88,wagneur91,gaubert98n,sturmfels03,CunBut04}. The general case has been less studied,
with the exception of the paper~\cite{helbig}, in which Helbig
established a max-plus analogue of Krein-Milman's 
theorem, showing that a non-empty 
compact convex subset of $(\R\cup\{-\infty\})^n$ is the closure 
of the convex hull of its set of extreme points, in the max-plus sense.

For conventional convex sets of finite dimension, however, a more precise 
result is true: the closure operator can be dispensed with,
since a classical theorem of Minkowski shows that a non-empty compact 
convex subset of a finite dimensional space is the 
convex hull of its set of extreme points. One may
ask whether the same is true for max-plus convex sets.
We show that the answer is positive, and establish
a max-plus analogue of Minkowski's theorem.

Note that the classical proof of Minkowski's
theorem cannot be transposed to the max-plus case.
The classical approach exploits 
the facial structure of convex sets. Recall that a {\em face} of a 
convex set is by definition the intersection of the 
convex set with a supporting hyperplane. For a conventional
convex set, one can show 
that the extreme points of the faces are extreme points of the 
set, and use this observation to prove Minkowski's 
theorem, by induction on the dimension of the convex set. 
This does not work in the max-plus case, because
an extreme point of a face may not be an extreme point of the set,
as shown in Example~\ref{cex} below. 
Hence, it does not seem possible to use Helbig's approach
to derive the results of the present paper.

In fact, we give a direct proof
of a Minkowski type theorem for max-plus convex
cones (Theorem~\ref{TheoMain}), from which we deduce the max-plus
Minkowski theorem (Theorem~\ref{TheorMinkowski}), and its generalization to
the case of unbounded convex sets (Theorem~\ref{TheorDesc}).
Finally, we deduce as a special case
a slightly more precise version of the ``basis theorem''
of Moller~\cite{moller88} and Wagneur~\cite{wagneur91}
for finitely generated max-plus convex
cones, Corollary~\ref{cor-basis}.

Finally, we note that the main results of the present paper,
Theorems~\ref{TheoMain}--\ref{TheorDesc}, have been announced
(without proof) in the survey paper~\cite{relmics}.
\section{Preliminaries}
In this section, we give basic definitions and establish
some elementary lemmas. 
To bring to light the analogy with classical convex analysis,
we shall use the following notation.
We denote by $\rmax$ the max-plus semiring.
We denote by $a\oplus b:=a\vee b$ the max-plus addition,
and by $ab:=a+b$ the max-plus multiplication. 
We set $\zero:=-\infty$, $\unit:=0$.
The set of vectors of size $n$ over $\rmax$ is denoted
by $\rmax^n$. A vector consisting only of $\zero$ entries
is denoted by $\zero$. By scalar, we mean an element of $\rmax$.
If $u,v\in \rmax^n$ and $\lambda\in \rmax$,
we set $u\oplus v:=u\vee v$, and we denote by $\lambda u$ 
the vector with entries $\lambda+u_i$. Max-plus convex sets 
and cones have been defined in the introduction. In the sequel,
for brevity, the term ``convex'' used without precisions
shall always be understood in
the max-plus sense. By ``cone'', we shall always mean
a (max-plus) convex cone.
\begin{definition}
Let $A$ be a subset of $\rmax^n$.

The {\em convex hull} of $A$, denoted by $\co (A)$, 
is the set of all (finite) {\em convex combinations} of elements of $A$.
These can be written as $\oplus_{k\in K}\alpha_ku^k$,
where $K$ is a finite set, 
$\left\{u^k\right\}_{k\in K}$ is a family of elements of 
$A$ and $\left\{\alpha_k\right\}_{k\in K}$ are scalars 
that satisfy $\oplus_{k\in K} \alpha_k =\unit$. 

The {\em cone} generated by $A$, denoted by $\cone (A)$, 
is the set of all (finite) linear combinations of elements of $A$.
These can be written as $\oplus_{k\in K}\alpha_ku^k$,
where $K$ is a finite set, 
$\left\{u^k\right\}_{k\in K}$ is a family of elements of $A$ and 
$\left\{\alpha_k\right\}_{k\in K}$ are scalars.

The {\em recession cone} of $A$ at a point $v\in A$
is defined by:
\[
\Rec_v (A):=\set{u\in \rmax^n}{v\oplus \lambda u\in A 
\mbox{ for all } \lambda \in \rmax } \enspace .
\]
\end{definition}

\begin{example}\label{ejemplo1}
The recession cone of the convex set $A$ of Figure~\ref{figure1}
at any point $v\in A$ is $ \Rec_v (A)=\cone (\left\{ (0,1),(2,0)\right\} )$.
It is shown on the right 
hand side of Figure~\ref{figure2}, below.

\begin{figure}
\begin{center}
\input{convexset1}
\end{center}
\caption{An unbounded max-plus convex set and three segments in general position contained in it.}
\label{figure1}
\end{figure}
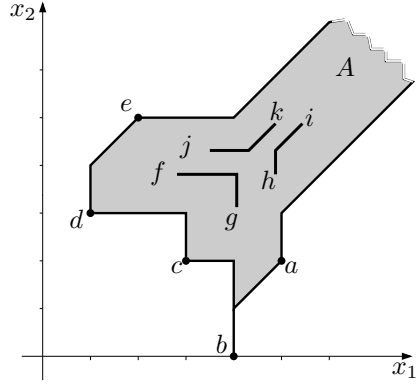

\end{example}

We equip $\rmax $ with the usual topology, which 
can defined by the metric: $(x,y) \rightarrow |e^x-e^y|$. 
The set $\rmax^n $ is equipped with the product topology. 
We denote by $\clo (A)$ the closure of a subset $A$ of $\rmax^n$.

The following lemmas give some properties of recession cones.

\begin{lemma}
Let $A$ be a closed subset of $\rmax^n$. Then the recession 
cone of $A$ at $v$ is closed for all $v\in A$.
\end{lemma} 
 
\begin{proof}
Let $v\in A$. For all $\lambda \in \rmax$, define the map 
$\varphi_\lambda:\rmax^n \rightarrow \rmax^n$ by 
$\varphi_\lambda(u)=v\oplus \lambda u$. Since $\varphi_\lambda$ 
is continuous,
$\Rec_v (A)=\cap_{\lambda \in \rmax} \varphi^{-1}_\lambda(A)$
is an intersection of closed sets, and so, it is closed.
\end{proof}

\begin{lemma}\label{prop1}
Let $A$ be a convex subset of $\rmax^n$ and $v,w\in A$. If 
$\beta v\leq w$ for some $\beta \neq\zero$, then 
$\Rec_v (A)\subset  \Rec_w (A)$. 
\end{lemma}

\begin{proof}
Let $u\in \Rec_v (A)$. We assume, 
without loss of generality, that $\beta \leq \unit$. Then,
for all $\lambda \in \rmax$, 
\[
w\oplus \lambda u = w\oplus \beta v\oplus \beta \beta^{-1}\lambda u = 
w\oplus \beta (v\oplus \beta^{-1} \lambda u)\in A
\]
since  $v\oplus \beta^{-1} \lambda u\in A$ and since $A$ is convex.
It follows that $u\in \Rec_w (A)$.
\end{proof}
\begin{propdef}\label{PropRecInde}
Let $A$ be a closed convex subset of $\rmax^n$. Then 
the recession cone of $A$ at $v$ is independent of $v\in A$. 
We denote it by $\Rec(A)$.
\end{propdef}
\begin{proof}
\typeout{Warning, support defined inside a proof}
Given $x\in \rmax^n$, we define the {\em support} of $x$ 
to be the set 
\[ \supp x:=\set{1\leq i\leq n}{x_i\neq \zero}
\enspace .\]
Observe first that if $u,v\in A$, then $u\oplus v\in A$
and the support of $u\oplus v$ is the union of the supports 
of $u$ and $v$. It follows that there is an element
$w\in A$ with maximum support, meaning that 
$\supp v \subset \supp w$ for all $v \in A$.  
Hence for every $v\in A$ there exists $\lambda_v \neq \zero$ 
such that $\lambda_v v \leq w$. Therefore, by Lemma~\ref{prop1}, it 
suffices to show that $\Rec_w(A) \subset \Rec_v(A)$ for every $v\in A$.

Let $\left\{\beta_r \right\}_{r\in \N} \subset \rmax$ 
be a sequence such that $\lim_{r\rightarrow \infty}\beta_r =\zero$ 
and $\zero< \beta_r\leq \unit $ for all $r\in \N$. 
If $u\in \Rec_w(A)$ and $\lambda \in \rmax$, then  
\[
v\oplus \lambda u =\lim_{r\rightarrow \infty}(v\oplus \beta_r (w\oplus \beta_r^{-1} \lambda u))
\] 
is a limit of elements of $A$ because 
$w\oplus \beta_r^{-1} \lambda u\in A$ for all 
$r\in \N$. Since $A$ is closed, it follows that 
$v\oplus \lambda u\in A$. Therefore, $u\in \Rec_v(A)$,
and so, $\Rec_w(A) \subset \Rec_v(A)$ for all $v\in A$. 
\end{proof}
\begin{remark} The closure assumption
in the previous proposition cannot be dispensed with.
Consider
$A=([\zero ,\unit ]\times \left\{ \zero \right\})\cup \R^2 =
\set{(x_1,\zero )}{x_1\leq \unit}\cup \R^2\subset \rmax^2$. 
Then,
\[
\Rec_{(\unit ,\zero )}(A)=\set{(u_1,u_2)\in \rmax^2}{u_2\neq \zero}\cup 
\left\{ (\zero ,\zero) \right\} 
\mbox { and }
\Rec_{(\unit ,\unit )}(A)=\rmax^2 \enspace .
\] 
\end{remark}
The following max-plus analogue of the notion of extreme point
was already used by Helbig~\cite{helbig}.
\begin{definition}[Extreme point]\label{DefExtPoint}
Let $A$ be a convex subset of $\rmax^n$. An element $x\in A$ is 
an {\em extreme point} of $A$ if for all $y,z\in A$ and 
$\alpha,\beta \in\rmax$ such that $\alpha \oplus \beta=\unit$, 
the following property is satisfied 
\begin{align}
x=\alpha y\oplus \beta z \implies x=y \mbox { or } x=z.
\label{e-ext}
\end{align}
The set of extreme points of $A$ will be denoted by $\ext (A)$.
\end{definition}
Thus, a point of $A$ is extreme if it cannot belong
to a segment of $A$ unless it is an end of this segment.
We warn the reader that due to the idempotency of addition,
the property~\eqref{e-ext}, with $\alpha\oplus \beta=\unit$
and $\alpha,\beta\neq \zero$ is not equivalent to
\[
x=\alpha y\oplus \beta z \implies x=y \mbox { and } x=z.
\]
\begin{remark}\label{PropTonta}
If $x\in A$ is an extreme point of $A$, then 
$ x=\alpha y\oplus \beta z$, with 
$\alpha \oplus \beta=\unit$ and $y,z\in A$, implies:
\[
(x=y,\alpha =\unit )\mbox{ or }(x=z,\beta =\unit ).
\]
Indeed, assume that $x=y$ but $\alpha <\unit$. Then, 
$\beta =\unit $. Assume by contradiction that $x\neq z$.
Then, we have $x_i>z_i$, for some $1\leq i\leq n$,
and so, $x_i=\alpha y_i \oplus z_i=\alpha x_i\oplus z_i<x_i$,
which is nonsense.
\end{remark}

When $C\subset \rmax^n$ is a cone, it is clear that its 
only extreme point is $\zero $. In this case, the relevant 
notion is that of extreme generator.

\begin{definition}[Extreme generator]
Let $C\subset \rmax^n$ be a cone. An element $x\in C$ is an 
{\em extreme generator} of $C$ if the following property is satisfied 
\[
x=y\oplus z,\;y,z\in C  \implies x=y \mbox { or } x=z. 
\] 
If $x$ is an extreme generator of $C$, 
then the set $\rmax x=\set{\lambda x}{\lambda \in \rmax}$ is 
an {\em extreme ray} of $C$. The set of extreme 
generators of $C$ will be denoted by $\eg (C)$.
\end{definition}

Extreme generators are called {\em join 
irreducible} elements in lattice theory.

\begin{remark}
It can be readily checked that every element of 
an extreme ray of $C$ is an extreme generator of $C$.
\end{remark}

\begin{example}\label{ejemplo2}
Let us consider the closed convex set 
$A\subset \rmax^2$ shown in Figure~\ref{figure1}. 
It can be easily seen that its extreme points are 
$a=(5,2)$, $b=(4,0)$, $c=(3,2)$, $d=(1,3)$ and 
$e=(2,5)$. The extreme rays of 
$\Rec (A)$ are $\rmax (0,1)$ and $\rmax (2,0)$ 
(see the right hand side of Figure~\ref{figure2}). 
\end{example}

The following construction will allow us to derive
results for convex sets as consequences
of results for cones.
\begin{lemma}\label{lemmaConvCone}
Let $A$ be a convex subset of $\rmax^n$. Then, the set
\[
\Gcone A =\set{(\lambda x, \lambda)}{x\in A,\lambda \in \rmax} 
\subset \rmax^{n+1}
\]
is a cone.
\end{lemma}

\begin{proof}
Let $\beta \in \rmax $ and 
$(\lambda_1 x^1,\lambda_1), (\lambda_2 x^2,\lambda_2)\in \Gcone A $, 
with $x^1,x^2\in A$ and $\lambda_1,\lambda_2\in \rmax $.
Assume, without loss of generality, that 
$\lambda:=\lambda_1\oplus \lambda_2\neq \zero$. 
Since $A$ is convex, 
$\lambda^{-1}(\lambda_1 x^1\oplus \lambda_2 x^2)\in A$,
and so, 
\[
(\lambda_1 x^1,\lambda_1)\oplus (\lambda_2 x^2,\lambda_2)=
(\lambda \lambda^{-1}
(\lambda_1 x^1\oplus \lambda_2 x^2),\lambda)
\in \Gcone A .
\]
Moreover, $C_A$ is obviously preserved by the multiplication
by a scalar.
\end{proof}
 
\begin{lemma}\label{lemmaClo}
If $A\subset \rmax^n$ is a convex set
then $\clo (A)$ is a convex set. The same is true
for cones.
\end{lemma}
\begin{proof}
This follows from the continuity of the functions
$(x,y)\rightarrow x\oplus y$ and $(\lambda ,x)\rightarrow \lambda x$.
\end{proof}
We next establish some properties of the cone $\Gcone A$.
\begin{proposition}\label{DesSigA}
If $A\subset \rmax^n$ is a closed convex set, then 
\[
\clo (\Gcone A )=
\Gcone A \cup (\Rec (A)\times \left\{ \zero \right\}).
\]
\end{proposition}

\begin{proof}
Let $(y,\alpha )\in \clo (\Gcone A )$.

Assume first that $\alpha \neq \zero$. Since 
$(y,\alpha )\in \clo (\Gcone A )$, there exists a sequence 
$\left\{(\lambda_r x^r,\lambda_r)\right\}_{r\in \N}\subset \Gcone A $ 
such that 
$\lim_{r\rightarrow \infty }(\lambda_r x^r,\lambda_r)=(y,\alpha )$. 
Then, as $\lim_{r\rightarrow \infty }\lambda_r=\alpha \neq \zero$
and $A$ is closed, we know that 
$x:=\lim_{r\rightarrow \infty }x^r=
\lim_{r\rightarrow \infty }\lambda_r^{-1}\lambda_rx^r=
\alpha^{-1} y$ belongs to $A$. Therefore, 
$(y,\alpha )=\lim_{r\rightarrow \infty }(\lambda_r x^r,\lambda_r)=
(\alpha x,\alpha )\in \Gcone A $. 

Assume now that $\alpha = \zero$. Let $x\in A$ and $\beta \in \rmax$. 
To prove that $(y,\alpha )\in (\Rec (A)\times \left\{ \zero \right\})$ 
it suffices to show that $x\oplus \beta y\in A$. As $x\in A$ we know 
that $(x,\unit )\in \Gcone A $. Using the fact that 
$\clo (\Gcone A )$ is a cone (by Lemmas~\ref{lemmaConvCone} 
and~\ref{lemmaClo}), it follows that 
$(x\oplus \beta y,\unit )=(x,\unit )\oplus \beta (y,\zero )
\in \clo (\Gcone A )$. Then, there exists a sequence 
$\left\{(\lambda_r x^r,\lambda_r)\right\}_{r\in \N}\subset \Gcone A $ 
such that 
$\lim_{r\rightarrow \infty }(\lambda_r x^r,\lambda_r)=(x\oplus \beta y,\unit )$. 
Therefore, 
\[
x\oplus \beta y=\lim_{r\rightarrow \infty }\lambda_r x^r =
(\lim_{r\rightarrow \infty }\lambda_r)(\lim_{r\rightarrow \infty } x^r)=
\lim_{r\rightarrow \infty } x^r \in \clo (A)=A. 
\]

Thus, $\clo (\Gcone A ) \subset 
\Gcone A \cup (\Rec (A)\times \left\{ \zero \right\})$.

Obviously $\Gcone A \subset \clo (\Gcone A )$. 
Let now $(y,\zero )\in \Rec (A)\times \left\{ \zero \right\}$. 
Take any $x\in A$. We know that $x\oplus \lambda y\in A$ for all 
$\lambda \in \rmax$. Then, if 
$\left\{ \lambda_r\right\}_{r\in \N}\subset \rmax$ is a sequence 
such that $\lim_{r\rightarrow \infty }\lambda_r^{-1}=\zero $, 
it follows that 
$(y,\zero )= \lim_{r\rightarrow \infty }(\lambda_r^{-1}(x\oplus \lambda_r y),\lambda_r^{-1})$ 
and therefore $(y,\zero )\in \clo (\Gcone A )$ since 
$(\lambda_r^{-1}(x\oplus \lambda_r y),\lambda_r^{-1})\in \Gcone A $ 
for all $r\in \N$.

Thus, $\Gcone A \cup (\Rec (A)\times \left\{ \zero \right\})
\subset \clo (\Gcone A )$.
\end{proof}

\begin{corollary}\label{CompConvSet}
If $A\subset \rmax^n$ is a compact convex set, then 
$\Gcone A \subset \rmax^{n+1}$ is a closed cone.
\end{corollary}
\begin{proof}
If $A$ is a compact subset of $\rmax^n$, it must be bounded from
above, and so $\Rec(A)=\{\zero\}$. By Proposition~\ref{DesSigA},
$\clo(C_A)=C_A\cup \{\zero\}=C_A$, and so, $C_A$ is closed.
\end{proof}

\begin{lemma}\label{lemmaProp1}
Let $A$ be a closed convex subset of $\rmax^n$. Then, 
\[
\eg (\clo (\Gcone A ))\cap (\Rec (A)\times \left\{ \zero \right\} )
= \eg (\Rec (A))\times  \left\{ \zero \right\} .
\]
\end{lemma}

\begin{proof}
Let $(x,\zero )\in \eg (\clo (\Gcone A ))\cap 
(\Rec (A)\times \left\{ \zero \right\} )$. Then, 
in particular, $x\in \Rec (A)$. Assume that $x=y\oplus z$, 
with $y,z\in \Rec(A)$. As $\clo (\Gcone A )=
\Gcone A \cup (\Rec (A)\times \left\{ \zero \right\})$ 
by Proposition~\ref{DesSigA}, we know that 
$(y,\zero ),(z,\zero ) \in \clo (\Gcone A )$. Then 
$x=y$ or $x=z$, since 
$(x,\zero )=(y,\zero )\oplus (z,\zero )$ and 
$(x,\zero )\in \eg (\clo (\Gcone A ))$. Therefore, 
$x\in \eg (\Rec (A))$ and 
$(x,\zero )\in \eg (\Rec (A))\times  \left\{ \zero \right\}$.

Let now $(x,\zero )\in \eg (\Rec (A))\times  \left\{ \zero \right\}$. 
Then, obviously $(x,\zero )\in \Rec (A)\times  \left\{ \zero \right\}$.
Assume that $(x,\zero )=(x^1,\lambda_1 )\oplus (x^2,\lambda_2 )$, 
with $(x^1,\lambda_1 ),(x^2,\lambda_2 )\in \clo (\Gcone A )$. 
Then, $x=x^1\oplus x^2$ and $\lambda_1=\lambda_2=\zero $. Therefore, as 
$\clo (\Gcone A )=\Gcone A \cup (\Rec (A)\times \left\{ \zero \right\})$ 
by Proposition~\ref{DesSigA}, it follows that $x^1,x^2\in \Rec(A)$. 
Finally, $x=x^1$ or $x=x^2$ since $x\in \eg (\Rec (A))$. 
Thus, $(x,\zero )\in  \eg (\clo (\Gcone A ))$.
\end{proof}

\begin{lemma}\label{lemmaProp2}
Let $A$ be a convex subset of $\rmax^n$. Then, 
\[
\eg (\clo (\Gcone A ))\cap \Gcone A \subset \eg (\Gcone A ).
\]
\end{lemma}

\begin{proof}
Obvious since $\Gcone A \subset \clo (\Gcone A )$.
\end{proof}

The following proposition relates extreme points and extreme rays.

\begin{proposition}\label{prop-rel}
Let $C\subset \rmax^n$ be a cone, let $\gamma \neq \zero$ and 
let $\psi :\rmax^n \rightarrow \rmax$ be a max-plus linear form, 
meaning that $\psi (x)=\oplus_{i=1}^{n}a_ix_i$ for some $a\in \rmax^n$. 
Assume that $\psi (x)\neq\zero$ for all $x\in C\setminus\{\zero\}$,
and define the convex set:
\[
\Sigma :=\set{x\in C}{\psi (x)=\gamma} \enspace. 
\]
Then, 
\[
\ext (\Sigma )=\eg (C)\cap \Sigma \enspace .
\]
\end{proposition}

\begin{proof}
Let $x,y\in \Sigma $ and $\alpha ,\beta \in \rmax$ 
be such that $\alpha \oplus \beta =\unit $. Then, as 
\[
\psi (\alpha x\oplus \beta y)=
\alpha \psi (x)\oplus \beta \psi (y)=
\gamma \alpha \oplus \gamma \beta =
\gamma (\alpha \oplus \beta )= \gamma  
\]
and obviously $\alpha x\oplus \beta y\in C$, 
it follows that $\alpha x\oplus \beta y\in \Sigma$. 
Therefore, $\Sigma$ is convex.

Let $x\in \ext (\Sigma )$. Assume that $x=y\oplus z$, for 
some $y,z\in C\setminus\{\zero\}$.
Then, $\psi (y)=\gamma$ or $\psi (z)=\gamma$ 
since $\gamma =\psi (x)=\psi (y) \oplus \psi (z)$. Suppose, 
without loss of generality, that $\psi (y)=\gamma$. As 
$x=y\oplus \psi (z)\gamma^{-1} \gamma \psi (z)^{-1}z$,
where clearly $\gamma \psi (z)^{-1}z\in \Sigma $ and 
$\psi (z)\gamma^{-1}\leq \unit $, 
we know that 
\[
x=y \mbox{ or } x= \gamma \psi (z)^{-1}z.
\]
Since $x\neq y$ implies $\psi (z) \gamma^{-1}=\unit $
(see Remark~\ref{PropTonta}), it follows that 
$x=y$ or $x=z$. Then, $x\in \eg (C)\cap \Sigma$.

Let now $x\in \eg (C)\cap \Sigma$. Suppose that 
$x=\alpha y\oplus \beta z$, with $y,z\in \Sigma $
and $\alpha \oplus \beta =\unit $. Since 
$x\in \eg (C)$, we know that $x=\alpha y$ or 
$x=\beta z$. Assume, without loss of generality, 
that $x=\alpha y$. Then, 
$\gamma =\psi (x)=\psi (\alpha y)=\alpha \gamma $ 
implies that $\alpha =\unit$, and so $x=y$ . Therefore, 
$x\in \ext (\Sigma )$.
\end{proof}

Note that the condition of the previous proposition 
is satisfied, in particular, when 
$\psi (x)=\oplus_{i=1}^{n}a_ix_i$ for some $a\in \R^n$.

\begin{corollary}\label{ExtGExtP}
Let $A$ be a convex subset of $\rmax^n$. Then, 
\[
\eg (\Gcone A )\cap (A\times \left\{ \unit \right\})
=\ext (A)\times \left\{ \unit \right\} .
\]
\end{corollary}
\begin{proof}
Consider the max-plus linear form $\psi$ on $\rmax^{n+1}$
defined by $\psi(z,\lambda)=\lambda$, for all $z\in \rmax^n$ and $\lambda\in \rmax$, take $\gamma:=\unit$, and apply Proposition~\ref{prop-rel}
to the cone $C_A\subset \rmax^{n+1}$. 
We deduce that $\eg(C_A)\cap \Sigma=\ext(\Sigma)$,
where $\Sigma:=\set{(x\lambda,\lambda)\in C_A}{\lambda=\unit}
= A\times \{\unit\}$. Since $\ext(A\times \{\unit\})=\ext(A)\times \{\unit\}$,
the corollary is proved.
\end{proof}
Let us recall that a cone $C\subset \rmax^n$ is
{\em finitely generated} if there exists a finite subset 
$A\subset \rmax^n$ such that $C=\cone (A)$. 
\begin{lemma}\label{lem-closed}
A finitely generated cone of $\rmax^n$ is closed.
\end{lemma}
\begin{proof}
Let $A=\{u^1,\ldots,u^m\}$ and $C=\cone(A)$.
We assume, without loss of generality, that $u^k \neq \zero$ for all 
$1\leq k\leq m$. Let $\psi(x):=\bigoplus_{1\leq i\leq n}a_i x_i$ denote a
linear form, such that $a_i>\zero$ for all $1\leq i\leq n$.
Then, $\psi(u^k)\neq\zero$, for all $1\leq k\leq m$.
Let $\left\{ x^r=\oplus_{k=1}^m \lambda^r_ku^k\right\}_{r\in \N}$ 
be a sequence of elements of $C$ such that 
$\lim_{r\rightarrow \infty} x^r=x$ for some 
$x\in \rmax^n$. 

Since $\lambda^r_k\psi(u^k)\leq \psi(x^r)$, and since $\psi(u^k)\neq\zero$,
$\lambda^r_k$ is bounded as $r$ tends to infinity. Hence,
we can assume, without loss of generality,
that there exists $\lambda_k\in \rmax$ such that 
$\lim_{r\rightarrow \infty}\lambda^r_k=\lambda_k$ 
for all $k=1,\ldots ,m$ (taking subsequences if 
necessary). Then,
\[ 
\lim_{r\rightarrow \infty} x^r =
\lim_{r\rightarrow \infty} \Big(\bigoplus_{k=1}^m \lambda^r_ku^k\Big)=
\bigoplus_{k=1}^m \lambda_ku^k \in C.
\]
Therefore, $C$ is closed.
\end{proof}

\section{Representation of max-plus convex sets in terms of extreme points and extreme generators}\label{sec-rep}

Now we prove the main results of this paper.

\begin{theorem}\label{TheoMain}
Let $C\subset \rmax^n$ be a non-empty closed cone. 
Then, every element of $C$ is the sum of at most
$n$ extreme generators of $C$, and so,
\[
C=\cone (\eg (C)).
\]
\end{theorem}
\begin{proof}
Let $x\in C$. 
For each $i\in \left\{1,\ldots ,n\right\}$ 
define the set 
\[
S_i=\set{u\in C}{u\leq x,u_i=x_i}= 
C\cap \set{u\in \rmax^n}{u\leq x,u_i=x_i}.
\]
As $\set{u\in \rmax^n}{u\leq x,u_i=x_i}$ 
is compact and $C$ is closed, 
we know that $S_i$ is a compact subset of $\rmax^n$ 
which is non-empty because $x\in S_i$. Therefore, 
$S_i$ has a minimal element $u^i$. 

We claim that $u^i$ is an extreme generator of $C$. 
Assume that $u^i=y\oplus z$ for some $y,z\in C$. 
Then, $u^i_i=y_i$ or $u^i_i=z_i$. Let us assume, 
without loss of generality, that $u^i_i=y_i$. 
Therefore, $y\in S_i$ since 
$y\leq u^i\leq x$ and $y_i=u^i_i=x_i$. Hence, $u^i=y$ 
since $y\leq u^i$ and $u^i$ is a minimal element of $S_i$. 
Thus, $u^i$ is an extreme generator of $C$. 
It is clear that $x=\oplus_{i=1}^n u^i$, 
and so, $x\in \cone (\eg (C))$. 

We have shown that $C\subset \cone (\eg (C))$.
The other inclusion is trivial.
\end{proof}

\begin{theorem}[Max-Plus Minkowski Theorem]\label{TheorMinkowski}
Let $A$ be a non-empty compact convex subset of $\rmax^n$. 
Then, every element of $A$ is the convex combination
of at most $n+1$ extreme points of $A$,
and so,
\[
A=\co (\ext (A)).
\]
\end{theorem}

\begin{proof}
Let $x\in A$. Define the cone $\Gcone A =
\set{(\lambda z, \lambda)}{z\in A,\lambda \in \rmax} 
\subset \rmax^{n+1}$ as in Lemma~\ref{lemmaConvCone}. 
Then, by Corollary~\ref{CompConvSet}, $\Gcone A$ is a closed cone 
and thus 
\[
\Gcone A =\cone (\eg (\Gcone A )) 
\]
by Theorem~\ref{TheoMain}.

As $(x,\unit )\in \Gcone A$, by Theorem~\ref{TheoMain} 
we know that there exist $n+1$ extreme generators of $\Gcone A$, 
namely $(\lambda_1 u^1,\lambda_1),\ldots ,
(\lambda_{n+1} u^{n+1},\lambda_{n+1})$, 
such that 
\[
(x,\unit )=\bigoplus_{k=1}^{n+1}(\lambda_k u^k,\lambda_k). 
\]
Hence, 
 \[
x=\bigoplus_{k=1}^{n+1}\lambda_k u^k,
\mbox {where }\bigoplus_{k=1}^{n+1}\lambda_k=\unit.  
\]
By Corollary~\ref{ExtGExtP}, we know that 
$(u^k,\unit )\in \Gcone A$ is an extreme generator of 
$\Gcone A$ if, and only if, $u^k\in A$ is an extreme 
point of $A$. This shows that $x$ is the convex
combination of at most $n+1$ extreme points of $A$.
It follows that $A \subset \co (\ext (A))$. The other inclusion
is trivial.
\end{proof}

\begin{theorem}\label{TheorDesc}
Let $A\subset \rmax^n$ be a non-empty closed convex set. 
Then, every element of $A$ is the sum of the convex combination
of $p$ extreme points of $A$, and of $q$ extreme generators of $\Rec(A)$,
with $p+q\leq n+1$, and so:
\[
A=\co (\ext (A))\oplus \Rec (A).
\]
\end{theorem}
Here, we denote by $\oplus$ the max-plus analogue of the 
Minkowski sum of two subsets,
which is defined as the set of max-plus sums of a
vector from the first set and of a vector from the second one.
\begin{proof}
Let $x\in A$. Define the cone $\Gcone A =
\set{(\lambda z, \lambda)}{z\in A,\lambda \in \rmax} 
\subset \rmax^{n+1}$ as in Lemma~\ref{lemmaConvCone}. 
Then, by Lemma~\ref{lemmaClo}, 
$\clo (\Gcone A)$ is a closed cone and thus 
\[
\clo (\Gcone A) =\cone (\eg (\clo (\Gcone A ))) 
\]
by Theorem~\ref{TheoMain}.

By Proposition~\ref{DesSigA} we know that 
\[
\clo (\Gcone A )=
\Gcone A \cup (\Rec (A)\times \left\{ \zero \right\}), 
\]
and then Lemmas~\ref{lemmaProp1} and~\ref{lemmaProp2} imply: 
\begin{eqnarray*}
\eg (\clo (\Gcone A ))& = &
\left[ \eg (\clo (\Gcone A ))\cap \Gcone A \right]\cup 
\left[ \eg (\clo (\Gcone A ))\cap (\Rec (A)\times \left\{ \zero \right\}) \right] \\ 
& \subset & \eg (\Gcone A ) \cup (\eg (\Rec (A))\times \left\{ \zero \right\}).
\end{eqnarray*}

Now, as $(x,\unit )\in\Gcone A\subset \clo (\Gcone A )$, 
by Theorem~\ref{TheoMain} we know that there exist 
a finite number of elements of $\eg (\Gcone A )$, 
namely $(\lambda_k u^k,\lambda_k)$ with $1\leq k \leq p$, 
and a finite number of elements of 
$\eg (\Rec (A))\times \left\{ \zero \right\}$, 
namely $(y^h,\zero)$ with $1\leq h \leq q$, such that 
\[
(x,\unit )=
\Big( \bigoplus_{1\leq k\leq p} (\lambda_k u^k,\lambda_k)\Big)
\oplus \Big(\bigoplus_{1\leq h\leq q} (y^h,\zero )\Big) ,
\] 
with $p+q\leq n+1$.
Therefore, 
\[
x=\Big( \bigoplus_{1\leq k\leq p} \lambda_k u^k \Big)
\oplus \Big(\bigoplus_{1\leq h\leq q} y^h\Big), 
\mbox { where }\bigoplus_{1\leq k\leq p}\lambda_k=\unit 
\]
and $y^h\in \Rec (A)$ for all
$1\leq h\leq q$.
By Corollary~\ref{ExtGExtP} we know that 
$(u^k,\unit )\in \Gcone A$ is an extreme generator of 
$\Gcone A$ if, and only if, $u^k\in A$ is an extreme point of $A$.
This shows that $x$ is the sum of the convex combination
of $p$ extreme points of $A$ and
of $q$ extreme generators of $\Rec(A)$ with $p+q\leq n+1$.
Hence, $A \subset \co (\ext (A))\oplus \Rec (A)$. 
The other inclusion is trivial.
\end{proof}

As a corollary of Theorem~\ref{TheoMain} 
we get a precise version of the ``basis theorem'' for finitely
generated cones. The first results of this
kind were obtained by Moller~\cite{moller88} and Wagneur~\cite{wagneur91}.
Several variants of this result have appeared 
in~\cite{gaubert98n,sturmfels03,CunBut04}. 

\begin{corollary}[Basis theorem]\label{cor-basis}
Let $C\subset \rmax^n$ be a finitely generated cone 
and $A\subset C$. Then, $C=\cone (A)$ if, and only if, 
$A$ contains at least one nonzero element of each 
extreme ray of $C$. 
\end{corollary}
\begin{proof}
This follows readily from Theorem~\ref{TheoMain} and
from Lemma~\ref{lem-closed}.
\end{proof}

\begin{example}
As an illustration of Theorem~\ref{TheorDesc}, 
let us consider once again the closed convex 
set $A\subset \rmax^2$ depicted in Figure~\ref{figure1}. 
We have already seen (Examples~\ref{ejemplo1} 
and~\ref{ejemplo2}) that 
$\ext (A)=\left\{ a,b,c,d,e\right\} $ and 
$\Rec (A)=\cone \left\{(0,1),(2,0)\right\} $. Then, 
\[
A= \co \left\{ a,b,c,d,e\right\}\oplus 
\cone \left\{(0,1),(2,0)\right\} 
\]
by Theorem~\ref{TheorDesc}. The sets 
$\co (\ext (A))$ and $\Rec (A)$ are depicted 
in Figure~\ref{figure2}.

\begin{figure}
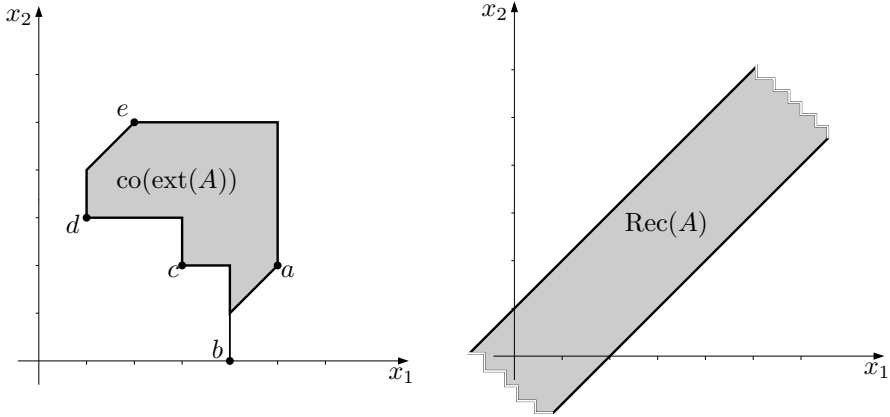

\begin{center}
\begin{tabular}[t]{cc}
\input convexset2 &
\input convexset3 
\end{tabular}
\end{center}
\caption{The sets $\co (\ext (A))$ and $\Rec (A)$ 
of Theorem~\ref{TheorDesc} for the unbounded
convex set depicted in Figure~\ref{figure1}.}
\label{figure2}
\end{figure}

\end{example}

\begin{example}\label{ex-notclosed}
The set of extreme points of a compact convex set may not be closed.
In the max-plus case, there are even counter-examples in 
dimension 2. Such a counter-example is shown in Figure~\ref{figure3},
where the set $A$ is given by 
\[
A = ([-2,0]\times \{0\})\cup (\{0\}\times [0,-2])
\cup \set{x\in \R^2}{-1\leq x_1+x_2,\; x_1,x_2\leq 0}
\enspace ,
\]
and 
\[
\ext(A)=
\{(-2,0),(0,-2)\}
\cup \set{x\in \R^2}{-1=x_1+x_2,\; x_1,x_2< 0}
\enspace .
\]
\end{example}
\begin{figure}
\begin{center}
\input{convexset4}
\end{center}
\caption{A convex subset of $\rmax^2 $ and its set of extreme points.}
\label{figure3}
\end{figure}
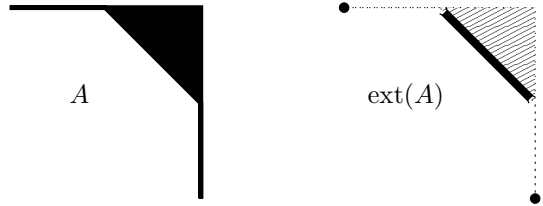
\begin{remark}
As in the classical case, the set of extreme points of a compact convex set
is a $\mathcal{G}_\delta$ set (a denumerable intersection of open sets).
Indeed, let $A$ be a non-empty compact convex subset of $\rmax^n$,
and let $d$ denote
any metric inducing the topology of $\rmax^n$. 
For all positive integers $k$, let 
\[ F_k:=\set{(x,y,z,\beta)\in A^3\times \rmax}{d(x,y)\geq 1/k,\; d(x,z)\geq 1/k,\; \beta\leq \unit,\; x=y\oplus \beta z}
\enspace .
\]
Let $\pi$ denote the projection sending $(x,y,z,\beta)$ to $x$.
Since $A$ is compact, $F_k$ is compact, and since $\pi$ is continuous,
$\pi(F_k)$ is compact. In particular, it is closed.
A point $x$ in $A$ is not extreme if and only
if there exist two points $y,z\in A$ that are both different from $x$,
and a scalar $\beta\leq \unit$, such that $x=y\oplus \beta z$.
The latter property means that $x$ belongs to some $\pi(F_k)$. So the
set of extreme points of $A$, which can be written as
$\cap_{k\geq 1} (\rmax^n\setminus \pi(F_k))$, is a $\mathcal{G}_\delta$
set.
\end{remark}

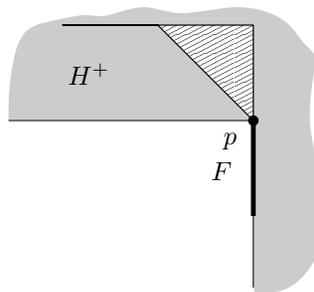
\begin{figure}
\begin{center}
\input{convexset4bis}
\end{center}
\caption{The point $p$ is an extreme point of the face $F$ but it is not an extreme point of the convex set}
\label{figure4bis}
\end{figure}

Let $a^+,a^-\in \rmax$ and let $\psi^+,\psi^-$ denote linear forms.
We call {\em half-space} a set of the form 
\[
H^+= \set{x\in \rmax^n}{\psi^+(x)\oplus a^+\geq \psi^-(x) \oplus a^-} \enspace.
\]
The {\em opposite half-space} $H^-$ is defined by reversing the inequality.
We say that $H^+$ is a {\em minimal supporting half-space} of $A$
if it contains $A$ and if it contains no other half-space containing
$A$. We define a {\em face} of a convex set $A$ to be the intersection
of $A$ with an half-space opposite to a minimal supporting half-space.
The following counter-example
shows that unlike in classical convex analysis,
the extreme points of faces are not necessarily
extreme points of the set.
\begin{example}\label{cex}
Consider the half-space
\[
H^+=\set{x\in\rmax^2}{x_1\oplus 1y_1 \geq 0} \enspace, 
\]
which is represented by the light gray region in Figure~\ref{figure4bis}.
One can check that this is a minimal supporting half-space of $A$
(see~\cite{cgqs03} or~\cite{joswig04} for a description
of max-plus half-spaces). Hence, $F:=A\cap H^-$ is a
face of $A$. This face is represented in bold on the figure.
The point $p=(0,-1)$ is an extreme point of $F$, but it is not
an extreme point of $A$.
\end{example}

\bibliographystyle{alpha} 
\bibliography{MPConvexbib} 
\end{document}

%% file: convexset1.tex
\begin{picture}(0,0)%
\includegraphics{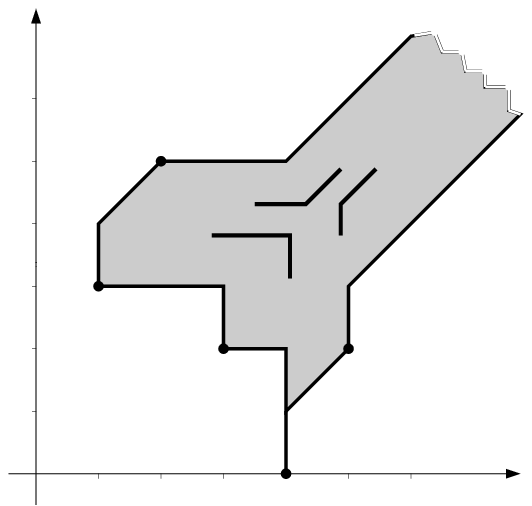}%
\end{picture}%
\setlength{\unitlength}{987sp}%
\begingroup\makeatletter\ifx\SetFigFont\undefined%
\gdef\SetFigFont#1#2#3#4#5{%
  \reset@font\fontsize{#1}{#2pt}%
  \fontfamily{#3}\fontseries{#4}\fontshape{#5}%
  \selectfont}%
\fi\endgroup%
\begin{picture}(10861,9569)(1576,-9383)
\put(3076,-5536){\makebox(0,0)[lb]{\smash{{\SetFigFont{10}{12.0}{\rmdefault}{\mddefault}{\updefault}{\color[rgb]{0,0,0}$d$}%
}}}}
\put(4351,-2611){\makebox(0,0)[lb]{\smash{{\SetFigFont{10}{12.0}{\rmdefault}{\mddefault}{\updefault}{\color[rgb]{0,0,0}$e$}%
}}}}
\put(5626,-6661){\makebox(0,0)[lb]{\smash{{\SetFigFont{10}{12.0}{\rmdefault}{\mddefault}{\updefault}{\color[rgb]{0,0,0}$c$}%
}}}}
\put(6751,-8686){\makebox(0,0)[lb]{\smash{{\SetFigFont{10}{12.0}{\rmdefault}{\mddefault}{\updefault}{\color[rgb]{0,0,0}$b$}%
}}}}
\put(8476,-6661){\makebox(0,0)[lb]{\smash{{\SetFigFont{10}{12.0}{\rmdefault}{\mddefault}{\updefault}{\color[rgb]{0,0,0}$a$}%
}}}}
\put(1576,-211){\makebox(0,0)[lb]{\smash{{\SetFigFont{10}{12.0}{\rmdefault}{\mddefault}{\updefault}{\color[rgb]{0,0,0}$x_2$}%
}}}}
\put(11176,-9211){\makebox(0,0)[lb]{\smash{{\SetFigFont{10}{12.0}{\rmdefault}{\mddefault}{\updefault}{\color[rgb]{0,0,0}$x_1$}%
}}}}
\put(9751,-1711){\makebox(0,0)[lb]{\smash{{\SetFigFont{10}{12.0}{\rmdefault}{\mddefault}{\updefault}{\color[rgb]{0,0,0}$A$}%
}}}}
\put(5851,-3661){\makebox(0,0)[lb]{\smash{{\SetFigFont{10}{12.0}{\rmdefault}{\mddefault}{\updefault}{\color[rgb]{0,0,0}$j$}%
}}}}
\put(8101,-2836){\makebox(0,0)[lb]{\smash{{\SetFigFont{10}{12.0}{\rmdefault}{\mddefault}{\updefault}{\color[rgb]{0,0,0}$k$}%
}}}}
\put(5101,-4336){\makebox(0,0)[lb]{\smash{{\SetFigFont{10}{12.0}{\rmdefault}{\mddefault}{\updefault}{\color[rgb]{0,0,0}$f$}%
}}}}
\put(6976,-5461){\makebox(0,0)[lb]{\smash{{\SetFigFont{10}{12.0}{\rmdefault}{\mddefault}{\updefault}{\color[rgb]{0,0,0}$g$}%
}}}}
\put(9001,-2986){\makebox(0,0)[lb]{\smash{{\SetFigFont{10}{12.0}{\rmdefault}{\mddefault}{\updefault}{\color[rgb]{0,0,0}$i$}%
}}}}
\put(7876,-4561){\makebox(0,0)[lb]{\smash{{\SetFigFont{10}{12.0}{\rmdefault}{\mddefault}{\updefault}{\color[rgb]{0,0,0}$h$}%
}}}}
\end{picture}%

%% file: convexset2.tex
\begin{picture}(0,0)%
\includegraphics{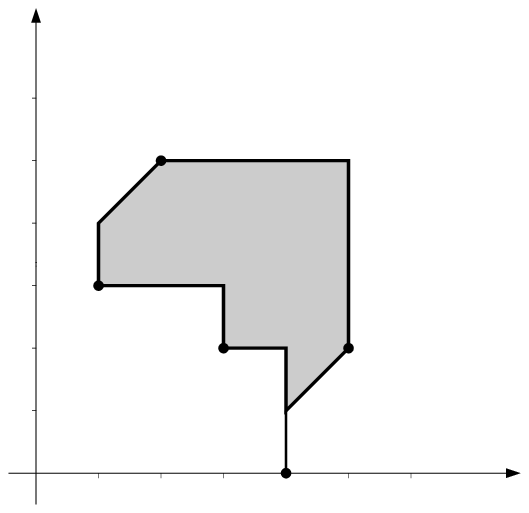}%
\end{picture}%
\setlength{\unitlength}{987sp}%
\begingroup\makeatletter\ifx\SetFigFont\undefined%
\gdef\SetFigFont#1#2#3#4#5{%
  \reset@font\fontsize{#1}{#2pt}%
  \fontfamily{#3}\fontseries{#4}\fontshape{#5}%
  \selectfont}%
\fi\endgroup%
\begin{picture}(10861,10309)(1576,-10123)
\put(3076,-5536){\makebox(0,0)[lb]{\smash{{\SetFigFont{10}{12.0}{\rmdefault}{\mddefault}{\updefault}{\color[rgb]{0,0,0}$d$}%
}}}}
\put(1576,-211){\makebox(0,0)[lb]{\smash{{\SetFigFont{10}{12.0}{\rmdefault}{\mddefault}{\updefault}{\color[rgb]{0,0,0}$x_2$}%
}}}}
\put(5626,-6661){\makebox(0,0)[lb]{\smash{{\SetFigFont{10}{12.0}{\rmdefault}{\mddefault}{\updefault}{\color[rgb]{0,0,0}$c$}%
}}}}
\put(6751,-8686){\makebox(0,0)[lb]{\smash{{\SetFigFont{10}{12.0}{\rmdefault}{\mddefault}{\updefault}{\color[rgb]{0,0,0}$b$}%
}}}}
\put(8476,-6661){\makebox(0,0)[lb]{\smash{{\SetFigFont{10}{12.0}{\rmdefault}{\mddefault}{\updefault}{\color[rgb]{0,0,0}$a$}%
}}}}
\put(11176,-9211){\makebox(0,0)[lb]{\smash{{\SetFigFont{10}{12.0}{\rmdefault}{\mddefault}{\updefault}{\color[rgb]{0,0,0}$x_1$}%
}}}}
\put(4351,-4411){\makebox(0,0)[lb]{\smash{{\SetFigFont{10}{12.0}{\rmdefault}{\mddefault}{\updefault}{\color[rgb]{0,0,0}co$($ext$(A))$}%
}}}}
\put(4351,-2611){\makebox(0,0)[lb]{\smash{{\SetFigFont{10}{12.0}{\rmdefault}{\mddefault}{\updefault}{\color[rgb]{0,0,0}$e$}%
}}}}
\end{picture}%

%% file: convexset3.tex
\begin{picture}(0,0)%
\includegraphics{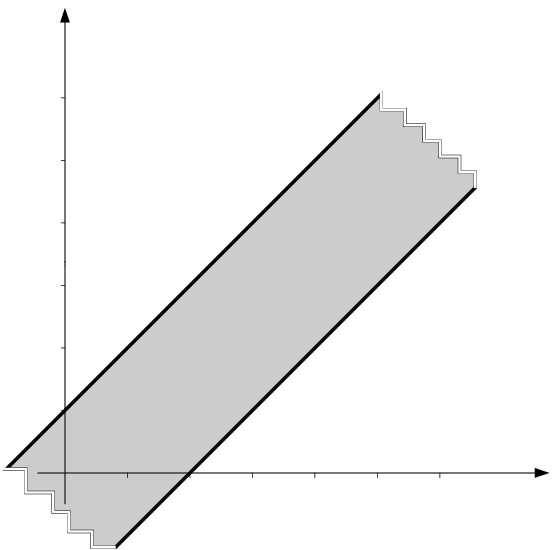}%
\end{picture}%
\setlength{\unitlength}{987sp}%
\begingroup\makeatletter\ifx\SetFigFont\undefined%
\gdef\SetFigFont#1#2#3#4#5{%
  \reset@font\fontsize{#1}{#2pt}%
  \fontfamily{#3}\fontseries{#4}\fontshape{#5}%
  \selectfont}%
\fi\endgroup%
\begin{picture}(11280,10426)(1157,-10240)
\put(11176,-9211){\makebox(0,0)[lb]{\smash{{\SetFigFont{10}{12.0}{\rmdefault}{\mddefault}{\updefault}{\color[rgb]{0,0,0}$x_1$}%
}}}}
\put(1576,-211){\makebox(0,0)[lb]{\smash{{\SetFigFont{10}{12.0}{\rmdefault}{\mddefault}{\updefault}{\color[rgb]{0,0,0}$x_2$}%
}}}}
\put(5176,-5611){\makebox(0,0)[lb]{\smash{{\SetFigFont{10}{12.0}{\rmdefault}{\mddefault}{\updefault}{\color[rgb]{0,0,0}Rec$(A)$}%
}}}}
\end{picture}%

%% file: convexset4.tex
\begin{picture}(0,0)%
\includegraphics{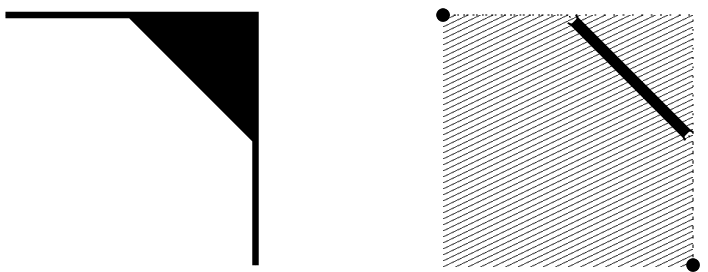}%
\end{picture}%
\setlength{\unitlength}{1973sp}%
\begingroup\makeatletter\ifx\SetFigFont\undefined%
\gdef\SetFigFont#1#2#3#4#5{%
  \reset@font\fontsize{#1}{#2pt}%
  \fontfamily{#3}\fontseries{#4}\fontshape{#5}%
  \selectfont}%
\fi\endgroup%
\begin{picture}(6722,2534)(1747,-4628)
\put(6301,-3361){\makebox(0,0)[lb]{\smash{{\SetFigFont{10}{12.0}{\rmdefault}{\mddefault}{\updefault}{\color[rgb]{0,0,0}ext$(A)$}%
}}}}
\put(2551,-3361){\makebox(0,0)[lb]{\smash{{\SetFigFont{10}{12.0}{\rmdefault}{\mddefault}{\updefault}{\color[rgb]{0,0,0}$A$}%
}}}}
\end{picture}%

%% file: convexset4bis.tex
\begin{picture}(0,0)%
\includegraphics{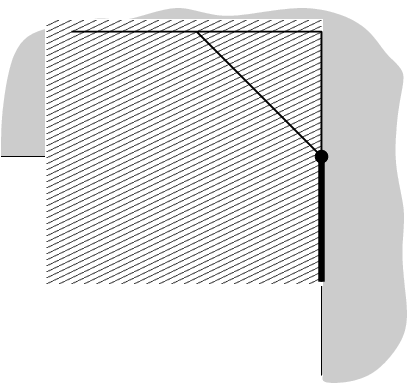}%
\end{picture}%
\setlength{\unitlength}{1973sp}%
\begingroup\makeatletter\ifx\SetFigFont\undefined%
\gdef\SetFigFont#1#2#3#4#5{%
  \reset@font\fontsize{#1}{#2pt}%
  \fontfamily{#3}\fontseries{#4}\fontshape{#5}%
  \selectfont}%
\fi\endgroup%
\begin{picture}(3907,3599)(5314,-5534)
\put(8026,-3661){\makebox(0,0)[lb]{\smash{{\SetFigFont{10}{12.0}{\rmdefault}{\mddefault}{\updefault}{\color[rgb]{0,0,0}$p$}%
}}}}
\put(7876,-4111){\makebox(0,0)[lb]{\smash{{\SetFigFont{10}{12.0}{\rmdefault}{\mddefault}{\updefault}{\color[rgb]{0,0,0}$F$}%
}}}}
\put(6076,-2911){\makebox(0,0)[lb]{\smash{{\SetFigFont{10}{12.0}{\rmdefault}{\mddefault}{\updefault}{\color[rgb]{0,0,0}$H^+$}%
}}}}
\end{picture}%